\documentclass[12pt,a4paper,reqno]{amsart} 
\usepackage{url}
\usepackage{fullpage} 

\newcommand{\eps}{\varepsilon}
\newcommand{\sameorder}{\asymp}
\newcommand{\dx}{\mathrm{d}} 
\newcommand{\I}{\mathcal{I}} 
\newcommand{\J}{\mathcal{J}} 
\newcommand{\C}{\mathbb{C}}

\newcommand{\R}{\mathbb{R}}
\newcommand{\Y}{\mathbb{Y}}

\newcommand{\Stilde}{\widetilde{S}}
\newcommand{\Odip}[2]{\mathcal{O}_{#1}\!\left(#2\right)\mathchoice{\!}{}{}{}}
\newcommand{\Odi}[1]{\Odip{}{#1}}

\allowdisplaybreaks

\renewcommand{\qedsymbol}{$\square$}

\newtheorem{Theorem}{Theorem}
\newtheorem{Lemma}{Lemma}

\begin{document} 

\title{A Ces\`aro Average of Hardy-Littlewood numbers}
%\date{\today}
\date{}
\author{Alessandro Languasco \& Alessandro Zaccagnini}
%
%    \subjclass is required.
\begin{abstract}
Let $\Lambda$ be the von Mangoldt function and
\(
  r_{\textit{HL}}(n)  =  \sum_{m_1 + m_2^2 = n} \Lambda(m_1),
\)
be the counting function for the Hardy-Littlewood numbers. Let $N$
be a sufficiently large integer. 
We prove that
\begin{align*}
  \sum_{n \le N} r_{\textit{HL}}(n) \frac{(1 - n/N)^k}{\Gamma(k + 1)}
  &=
  \frac{\pi^{1 / 2}}2 \frac{N^{3 / 2}}{\Gamma(k + 5 / 2)}
  -
  \frac 12 \frac{N}{\Gamma(k + 2)}
  -
  \frac{\pi^{1 / 2}}2
  \sum_{\rho}
    \frac{\Gamma(\rho)}{\Gamma(k + 3 / 2 + \rho)} N^{1 / 2 + \rho} \\
  &+
  \frac 12
  \sum_{\rho}
    \frac{\Gamma(\rho)}{\Gamma(k + 1 + \rho)} N^{\rho}
  +
  \frac{N^{3 / 4 - k / 2}}{\pi^{k + 1}}
  \sum_{\ell \ge 1}
    \frac{J_{k + 3 / 2} (2 \pi \ell N^{1 / 2})}{\ell^{k + 3 / 2}} \\
  &-
  \frac{N^{1 / 4 - k / 2}}{\pi^k}
  \sum_{\rho} \Gamma(\rho) \frac{N^{\rho / 2}}{\pi^\rho}
    \sum_{\ell \ge 1}
      \frac{J_{k + 1 / 2 + \rho} (2 \pi \ell N^{1 / 2})}
           {\ell^{k + 1 / 2 + \rho}}
  +
  \Odip{k}{1}.
\end{align*}
for $k > 1$, where $\rho$ runs over
the non-trivial zeros of the Riemann zeta-function $\zeta(s)$
and $J_{\nu} (u)$ denotes the Bessel function of complex order $\nu$
and real argument $u$.
\end{abstract}
\maketitle

\section{Introduction}

We continue our recent work on additive problems with prime summands. In \cite{LanguascoZ2012a} we studied the \emph{average} number of
representations of an integer as a sum of two primes, whereas in
\cite{LanguascoZ2012b} we considered individual integers.
In this paper we study a Ces\`aro weighted \emph{explicit} formula for
Hardy-Littlewood numbers (integers that can be written as
a sum of a prime and a square) and the goal is similar to the one in
\cite{LanguascoZ2012f}, that is, we want to obtain an asymptotic 
formula with the expected main
term and one or more terms that depend explicitly on the zeros of the
Riemann zeta-function.
Letting
\begin{equation}
\label{r2-def}
  r_{\textit{HL}}(n)  =  \sum_{m_1 + m_2^2 = n} \Lambda(m_1),
\end{equation}
the main result of the paper is the following theorem.

\begin{Theorem}
\label{CesaroHL-average} 
Let $N$ be a sufficiently large integer.
We have
\begin{align*}
   \sum_{n \le N} r_{\textit{HL}}(n) \frac{(1 - n/N)^k}{\Gamma(k + 1)}
  &=
  \frac{\pi^{1 / 2}}2 \frac{N^{3 / 2}}{\Gamma(k + 5 / 2)}
  -
  \frac 12 \frac{N}{\Gamma(k + 2)}
  -
  \frac{\pi^{1 / 2}}2
  \sum_{\rho}
    \frac{\Gamma(\rho)}{\Gamma(k + 3 / 2 + \rho)} N^{1 / 2 + \rho} \\
  &+
  \frac 12
  \sum_{\rho}
    \frac{\Gamma(\rho)}{\Gamma(k + 1 + \rho)} N^{\rho}
  +
\frac{N^{3 / 4 - k / 2}}{\pi^{k + 1}}
  \sum_{\ell \ge 1}
    \frac{J_{k + 3 / 2} (2 \pi \ell N^{1 / 2})}{\ell^{k + 3 / 2}} \\
  &-
  \frac{N^{1 / 4 - k / 2}}{\pi^k}
  \sum_{\rho} \Gamma(\rho) \frac{N^{\rho / 2}}{\pi^\rho}
    \sum_{\ell \ge 1}
      \frac{J_{k + 1 / 2 + \rho} (2 \pi \ell N^{1 / 2})}
           {\ell^{k + 1 / 2 + \rho}}
  +
  \Odip{k}{1}.
\end{align*}
for $k > 1$, where $\rho$ runs over
the non-trivial zeros of the Riemann zeta-function $\zeta(s)$
and $J_{\nu} (u)$ denotes the Bessel function of complex order $\nu$
and real argument $u$.
\end{Theorem}

Similar averages of arithmetical functions are common in the 
literature, see, e.g., Chan\-dra\-sekharan-Narasimhan
\cite{ChandrasekharanN1961} and Berndt \cite{Berndt1975} 
who built on earlier classical works (Hardy, Landau, Walfisz and others). 
In their setting the 
generalized Dirichlet series associated to the arithmetical 
function satisfies a suitable functional equation and this leads 
to an asymptotic formula containing Bessel functions of 
real order and argument.
In our case we have no functional equation, and, as far as we know, 
it is the first time that Bessel functions with complex order arise 
in a similar problem. Moreover, from a technical point of
view, the estimates of such Bessel functions are harder to perform 
than the ones already present in the Number Theory literature since 
the real argument and the complex order are both unbounded while, in previous 
papers, either the real order or the argument is bounded.

The method we will use in this additive problem is based on a formula due to Laplace
\cite{Laplace1812}, namely
\begin{equation}
\label{Laplace-transf}
  \frac 1{2 \pi i}
  \int_{(a)} v^{-s} e^v \, \dx v
  =
  \frac1{\Gamma(s)},
\end{equation}
where $\Re(s) > 0$ and $a > 0$, see, e.g., formula 5.4(1) on page 238
of \cite{ErdelyiMOT1954a}.
%\footnote{KP lo chiamano Bateman's project}
In the following we will need the general case of \eqref{Laplace-transf}
which can be found in de Azevedo Pribitkin \cite{Azevedo2002}, 
formulae (8) and (9):
\begin{equation}
\label{Laplace-eq-1}
  \frac1{2 \pi}
  \int_{\R} \frac{e^{i D u}}{(a + i u)^s} \, \dx u
  =
  \begin{cases}
    \dfrac{D^{s - 1} e^{- a D}}{\Gamma(s)}
    & \text{if $D > 0$,} \\
    0
    & \text{if $D < 0$,}
  \end{cases}
\end{equation}
which is valid for $\sigma = \Re(s) > 0$ and $a \in \C$ with
$\Re(a) > 0$, and
\begin{equation}
\label{Laplace-eq-2}
  \frac1{2 \pi}
  \int_{\R} \frac 1{(a + i u)^s} \, \dx u
  =
  \begin{cases}
    0     & \text{if $\Re(s) > 1$,} \\
    1 / 2 & \text{if $s = 1$,}
  \end{cases}
\end{equation}
for $a \in \C$ with $\Re(a) > 0$.
Formulae \eqref{Laplace-eq-1}-\eqref{Laplace-eq-2} enable us to write
averages of arithmetical functions by means of line integrals as we
will see in \S\ref{settings} below.

We will also need  Bessel functions of complex order $\nu$ 
and real argument $u$.
For their definition and main properties we
refer to Watson \cite{Watson1966}.
In particular, equation (8) on page 177 gives the
Sonine representation:
\begin{equation}
  \label{Bessel-def}
  J_\nu(u)
  :=
  \frac{(u / 2)^\nu}{2 \pi i}
  \int_{(a)} s^{- \nu - 1} e^s e^{- u^2 / 4 s} \, \dx s,
\end{equation}
where $a > 0$ and $u,\nu \in \C$ with $\Re(\nu) > -1$.
We will use also 
a Poisson integral formula (see 
%eq.~7.12.(8) on page  81 of Erdelyi-Magnus-Oberhettinger-Tricomi \cite
%{ErdelyiMOT1954b} or  eq.~(5.10.3) of Lebedev \cite{Lebedev1972} or  
eq.~(3) on page 48 of  
%Watson 
\cite{Watson1966}),  i.e.,
\begin{equation}
\label{Poisson-int-rep}
J_\nu(u)
: =
\frac{2(u/2)^{\nu}}{\pi^{1/2}\Gamma(\nu+1/2)}
\int_{0}^{1} (1-t^2)^{\nu-1/2} \cos (ut)\ \dx t
\end{equation}
which holds for $\Re (\nu) > -1/2$ and $u\in \C$.
An asymptotic estimate we will need is  
\begin{equation}
\label{Lebedev-asymp}
J_\nu(u)
=
\Bigl(\frac{2}{\pi u}\Bigr)^{1/2}
\cos \Bigl(u -\frac{\pi \nu}{2} -\frac{\pi}{4}\Bigr)
+
\Odip{\vert \nu \vert}{u^{-5/2}}
\end{equation}
which follows from 
%eq.~(5.11.6) of  Lebedev \cite{Lebedev1972} or 
eq.~(1) on page 199 of Watson \cite{Watson1966}.
 
As in \cite{LanguascoZ2012f}, we combine this approach with line
integrals with the classical methods dealing with infinite sums over
primes, exploited by Hardy \& Littlewood (see \cite{HardyL1916} and
\cite{HardyL1923}) and by Linnik \cite{Linnik1946}.
The main difference here is that the problem naturally involves the 
modular relation for the complex theta function,
see eq.~\eqref{func-eq-theta}; the presence of the Bessel functions
in our statement strictly depends on such modularity relation. 
It is worth mentioning that it is not clear how to get such ``modular'' terms
using the finite sums approach for the Hardy-Littlewood function 
$r_{\textit{HL}}(n)$.

We thank A.~Perelli and J.~Pintz for several conversations on this topic.

\section{Settings}
\label{settings}

We need $k> 0$ in this section.
Let  $z = a + i y$ with $a > 0$,
\begin{equation}
\label{Stilde-omega-def}
  \Stilde(z)
  =
  \sum_{m \ge 1} \Lambda(m) e^{- m z}
 \quad
 \text{and}
  \quad
  \omega_2(z)
  %&
  =
  \sum_{m \ge 1} e^{-m^2 z}.
\end{equation}
Letting further  
$\theta(z)  =   \sum_{m = -\infty}^{+\infty} e^{-m^2 z}$,
we notice that $\theta(z) = 1 + 2 \omega_2(z)$ and, recalling the
functional equation for $\theta$ (see, e.g., Proposition VI.4.3 of 
Freitag-Busam \cite[page 340]{FreitagB2009}):
\begin{equation}
\label{func-eq-theta}
  \theta(z)
  =
  \Bigl( \frac \pi z \Bigr)^{1/2}
  \theta\Bigl( \frac{\pi^2} z \Bigr),
\end{equation}
we immediately get 
\begin{equation}
\label{func-eq-omega}
  \omega_2(z)
  =
  \frac 12
  \Bigl( \frac \pi z \Bigr)^{1 / 2}
  -
  \frac12
  +
  \Bigl( \frac \pi z \Bigr)^{1 / 2}
  \omega_2 \Bigl( \frac {\pi^2} z \Bigr).
\end{equation}
Recalling \eqref{r2-def}, we can write
\[
  \Stilde(z) \omega_2(z)
  =
  \sum_{m_1 \ge 1} \sum_{m_2 \ge 1} \Lambda(m_1) e^{-(m_1 + m_2^2) z}
  =
  \sum_{n \ge 1} r_{\textit{HL}}(n) e^{- n z}
\]
and, by   \eqref{Laplace-eq-1}-\eqref{Laplace-eq-2}, we see  that
\begin{equation}
\label{first-step}
\sum_{n \le N} r_{\textit{HL}}(n) \frac{(N - n)^k}{\Gamma(k + 1)}
 =
\sum_{n \ge 1} r_{\textit{HL}}(n) 
\Bigl(
 \frac{1}{2 \pi i}
  \int_{(a)} e^{(N- n)z} z^{- k - 1}  \, \dx z
  \Bigr).
\end{equation}

Our first goal is to exchange the series with the line integral in
\eqref{first-step}. To do so
we have to recall that the Prime Number Theorem (PNT) is equivalent, via
Lemma~\ref{Linnik-lemma2} below, to the statement
\begin{equation*}
%\label{PNT-equiv}
  \Stilde(a)
  \sim
  a^{-1}
  \qquad\text{for $a \to 0+$,}
\end{equation*}
which is classical: for the proof see for instance Lemma~9 in Hardy \&
Littlewood \cite{HardyL1923}.
We will also use the inequality
\begin{equation}
\label{omega-estim}
  \vert \omega_2(z)\vert 
  \le
  \omega_2(a)  
  \le
  \int_{0}^{\infty} e^{-at^{2}} \dx t
   \le
   a^{- 1 / 2}
  \int_{0}^{\infty} e^{-v^{2}} \dx v
  \ll
  a^{- 1 / 2}
\end{equation}
from which we immediately get
\begin{align*}
  \sum_{n \ge 1} \bigl\vert  r_{\textit{HL}}(n) e^{- n z} \bigr\vert 
  &=
  \sum_{n \ge 2}   r_{\textit{HL}}(n) e^{- n a}  
=
  \Stilde(a) \omega_2(a)
  \ll
  a^{- 3 / 2}.
\end{align*}
Taking into account the estimates
\begin{equation}
\label{z^-1}
  \vert z \vert^{-1}
  \sameorder
  \begin{cases}
    a^{-1}   &\text{if $\vert y \vert \le a$,} \\
    \vert y \vert^{-1} &\text{if $\vert y \vert \ge a$,}
  \end{cases}
\end{equation}
 where $f\sameorder g$ means $g \ll f \ll g$, and
\[
  \vert e^{N z} z^{- k - 1}\vert 
  \sameorder
   e^{N a}
  \begin{cases}
    a^{- k - 1} &\text{if $\vert y \vert \le a$,} \\
    \vert y \vert^{- k - 1} &\text{if $\vert y \vert \ge a$,}
  \end{cases}
\]
we have
\begin{align*}
  \int_{(a)} \vert e^{N z} z^{- k - 1}\vert  \,
    \vert 
     \Stilde(z) \omega_2(z)
    \vert  \, \vert \dx z \vert
  &\ll
  a^{- 3 / 2} e^{N a}
  \Bigl(
    \int_{-a}^a a^{- k - 1} \, \dx y
    +
    2
    \int_a^{+\infty} y^{- k - 1} \, \dx y
  \Bigr) \\
  &\ll
   a^{- 3 / 2} e^{N a}
  \Bigl( a^{-k} + \frac{a^{-k}}k \Bigr),
\end{align*}
but the last estimate is valid only if $k > 0$.
So, for $k > 0$, we can exchange the line integral with the sum over $n$
in \eqref{first-step} thus getting
\begin{equation}
\label{main-form-omega}
  \sum_{n \le N} r_{\textit{HL}}(n) \frac{(N - n)^k}{\Gamma(k + 1)}
  =
  \frac 1{2 \pi i}
  \int_{(a)} e^{N z} z^{- k - 1} \Stilde(z) \omega_2(z) \, \dx z.
\end{equation}
This is the fundamental relation for the method.

\section{Inserting zeros and modularity}

We need $k>  1/2$ in this section.
The treatment of the integral at the right hand side
of \eqref{main-form-omega} requires Lemma \ref{Linnik-lemma2}.
Letting $E(a,y)$ be the error term in \eqref{expl-form-err-term-strong},
formula \eqref{main-form-omega} becomes
\begin{align*}
  \sum_{n \le N} r_{\textit{HL}}(n) \frac{(N - n)^k}{\Gamma(k + 1)}
  &=
  \frac 1{2 \pi i}
  \int_{(a)}
    \Bigl( \frac 1z - \sum_{\rho} z^{-\rho} \Gamma(\rho) \Bigr)
    \omega_2(z) e^{N z} z^{- k - 1} \, \dx z \\
  &\qquad+
  \Odi{\int_{(a)}
    \vert E(a,y)\vert  \, \vert e^{N z}\vert  \, \vert z\vert ^{- k - 1} \vert \omega_2(z)\vert  \, \vert \dx z\vert}.
\end{align*}
Using \eqref{omega-estim}-\eqref{z^-1}
and \eqref{expl-form-err-term-strong} we see that the error term is
\begin{align*}
%&\ll
% a^{- 1 / 2} e^{Na} 
% \Bigl( 
%  \int_{-a}^{a} \vert z \vert^{{-k-1/2}} \dx y
%+
% \int_{a}^{+\infty} \vert z \vert^{{-k-1/2}}\log^2(y/a) \, \dx y 
% \Bigr)
% \\
&\ll
 a^{- 1 / 2} e^{Na} 
 \Bigl( 
 \int_{-a}^{a} a^{{-k-1/2}} \dx y 
 +
 \int_{a}^{+\infty} y^{{-k-1/2}} \log^2(y/a)\,  \dx y 
 \Bigr)
 \\
&\ll_k
  e^{N a} a^{-k}
  \Bigl( 1+
 \int_{1}^{+\infty}v^{{-k-1/2}} \log^2 v \, \dx v
   \Bigr)
\ll_k
  e^{N a} a^{-k},
\end{align*}
provided that $k>1/2$.
Choosing $a = 1 / N$, the previous estimate becomes
\(
  \ll_k
  N^{k}.
\)
Summing up, for  $k>  1/2$, we can write
\begin{equation}
\label{first-step-hl}
  \sum_{n \le N} r_{\textit{HL}}(n) \frac{(N - n)^k}{\Gamma(k + 1)}
  =
  \frac 1{2 \pi i}
  \int_{(1 / N)}
    \Bigl(\frac 1z
          -
          \sum_{\rho} z^{-\rho} \Gamma(\rho)
    \Bigr)
    \omega_2(z) e^{N z} z^{- k - 1} \, \dx z
  +
  \Odip{k}{N^{k}}.
\end{equation}
We now insert \eqref{func-eq-omega} into \eqref{first-step-hl}, so that
the integral on the right-hand side of \eqref{first-step-hl} becomes
\begin{align}
\notag
  &
  \frac 1{2 \pi i}
  \int_{(1 / N)}
    \Bigl( \frac 1z - \sum_{\rho} z^{-\rho} \Gamma(\rho) \Bigr)
    \Bigl( \frac 12 \Bigl( \frac \pi z \Bigr)^{1/2} - \frac12 \Bigr)
    e^{N z} z^{- k - 1} \, \dx z \\
\notag
  &\qquad+
  \frac 1{2 \pi i}
  \int_{(1 / N)}
    \Bigl( \frac \pi z \Bigr)^{1/2}
    \Bigl( \frac 1z - \sum_{\rho} z^{-\rho} \Gamma(\rho) \Bigr)
    \omega_2 \Bigl( \frac{\pi^2}z \Bigr)
     e^{N z} z^{- k - 1} \, \dx z \\
\label{hl-splitting}
  &=
  \I_1 + \I_2,
\end{align}
say. We now proceed to evaluate $\I_1$ and $\I_2$.

\section{Evaluation of $\I_1$}
We need $k>  1/2$ in this section.
By a  direct computation we can write that
\begin{align*}
  \I_1
  &=
  \frac 1{4 \pi i}
  \int_{(1 / N)}
    \Bigl(  \frac{\pi^{1/2}}{z^{1/2}} - 1 \Bigr)
   e^{N z} z^{- k - 2}  \, \dx z
  -
  \frac {\pi^{1/2}}{4 \pi i} 
  \int_{(1 / N)}
    \sum_{\rho} \Gamma(\rho)  e^{N z} z^{- k -\rho - 3/2} \, \dx z \\
  &\qquad+
  \frac 1{4 \pi i} 
  \int_{(1 / N)}
    \sum_{\rho}  \Gamma(\rho) e^{N z}  z^{- k - \rho - 1}  \, \dx z 
  =
  \J_1 + \J_2 + \J_3,
\end{align*}
say. We see now how to evaluate $\J_1$, $\J_2$ and $\J_3$.

\subsection{Evaluation of $\J_1$}

Using the substitution $s=Nz$, 
by \eqref{Laplace-transf} we immediately have
\begin{align}
\notag
  \J_1 
  &=
  \frac{\pi^{1/2}}2 N^{k + 3/2}
  \frac 1{2 \pi i}
  \int_{(1)} e^s s^{- k - 5 / 2} \, \dx s
  -
  \frac 12 N^{k + 1}
  \frac 1{2 \pi i}
  \int_{(1)} e^s s^{- k - 2} \, \dx s \\
  \label{J1-eval}
  &=
  \frac{\pi^{1/2}}2 \frac{N^{k + 3/2}}{\Gamma(k + 5/2)}
  -
  \frac 12 \frac{N^{k + 1}}{\Gamma(k + 2)}.
\end{align}

\subsection{Evaluation of $\J_2$}

Exchanging the sum over $\rho$ with the integral (this can be done for $k>0$,
see  \S\ref{exchange-rho-integral}) and using the substitution $s=Nz$, we have
\begin{align}
\notag
  \J_2
  &=
  -
  \frac{\pi^{1/2}}2
  \sum_{\rho} \Gamma(\rho)
    \frac 1{2 \pi i}
    \int_{(1 / N)} e^{N z}  z^{-k -\rho - 3/2} \, \dx z \\
    \notag
  &=
  -
  \frac{\pi^{1/2}}2
  \sum_{\rho} \Gamma(\rho) N^{k + \rho + 1/2}
    \frac 1{2 \pi i}
    \int_{(1)} e^s s^{-k -\rho - 3/2}  \, \dx s   \\
    \label{J2-eval}
  &=
  -
  \frac{\pi^{1/2}}2
  \sum_{\rho}
    \frac{\Gamma(\rho)}{\Gamma(k + 3 / 2 + \rho)} N^{k + 1 / 2 + \rho},
\end{align}
again by \eqref{Laplace-transf}.
By the Stirling formula \eqref{Stirling}, 
we remark that the series in $\J_2$ converges absolutely 
 for  $k>-1/2$. 

\subsection{Evaluation of $\J_3$}

Arguing as in \S\ref{exchange-rho-integral} with $-k-1$ which plays 
the role of $-k-3/2$ there, we see that we can exchange the sum with the integral 
provided that $k>1/2$. Hence, performing again the usual substitution $s=Nz$,
we can write
\begin{equation}
\label{J3-eval}
  \J_3 
  =
  \frac 12
  \sum_{\rho} \Gamma(\rho) N^{k + \rho}
    \frac 1{2 \pi i}
    \int_{(1)} e^s s^{- k - 1 - \rho} \, \dx s \\
  =
  \frac 12
  \sum_{\rho}
     \frac{\Gamma(\rho)}{\Gamma(k + 1 + \rho)} N^{k + \rho}.
\end{equation}
By the Stirling formula \eqref{Stirling}, 
we remark that the series in $\J_3$ converges absolutely 
 for  $k>0$.

\medskip
Summing up, by \eqref{J1-eval}-\eqref{J3-eval} and for $k>  1/2$ we get
\begin{align}
\notag
  \I_1
  &=
  \frac{\pi^{1 / 2}}2 \frac{N^{k + 3 / 2}}{\Gamma(k + 5 / 2)}
  -
  \frac 12 \frac{N^{k + 1}}{\Gamma(k + 2)}
  -
  \frac{\pi^{1 / 2}}2
  \sum_{\rho}
    \frac{\Gamma(\rho)}{\Gamma(k + 3 / 2 + \rho)} N^{k + 1 / 2 + \rho} \\
  &\qquad+
\label{final-eval-I1}
  \frac 12
  \sum_{\rho}
    \frac{\Gamma(\rho)}{\Gamma(k + 1 + \rho)} N^{k + \rho}.
\end{align}

\section{Evaluation of $\I_2$ and conclusion of the proof of Theorem \ref{CesaroHL-average}}
 
 We need $k>  1$ in this section.
Using the definition of $\omega_2 (\pi^2/ z )$, 
see \eqref{Stilde-omega-def}, we have
\begin{align}
\notag
  \I_2
  &=
  \frac 1{2 \pi i}
  \int_{(1 / N)}
    \Bigl( \frac \pi z \Bigr)^{1/2}
    \Bigl( \sum_{\ell \ge 1} e^{- \ell^2 \pi^2 / z} \Bigr)
     e^{N z} z^{- k - 2} \, \dx z \\
  &
  \label{I2-splitting}
  \qquad
  - \frac 1{2 \pi i}
  \int_{(1 / N)}
    \Bigl( \frac \pi z \Bigr)^{1/2}
    \Bigl( \sum_{\ell \ge 1} e^{- \ell^2 \pi^2 / z} \Bigr)
    \Bigl( \sum_{\rho} z^{-\rho} \Gamma(\rho) \Bigr)
     e^{N z} z^{- k - 1} \, \dx z 
  =
  \J_4
  +
  \J_5,
\end{align}
say. We see now how to evaluate  $\J_4$ and $\J_5$.

\subsection{Evaluation of $\J_4$}

By means of the substitution $s = N z$, since
the exchange  is justified   in \S\ref{exchange-ell-integral} for $k> -1/2$, we get
\begin{equation*}
%\label{J4-def}
  \J_4 
  =
  \pi^{1 / 2}
  N^{k + 3 / 2}
  \sum_{\ell \ge 1}
    \frac 1{2 \pi i}
    \int_{(1)}
      e^s e^{- \ell^2 \pi^2 N / s} s^{- k - 5 / 2} \, \dx s.
\end{equation*}
Setting $u = 2 \pi \ell N^{1/2}$ in \eqref{Bessel-def} we obtain
\begin{equation}
\label{J-nu}
  J_\nu \bigl( 2 \pi \ell N^{1/2} \bigr)
  =
  \frac{(\pi \ell N^{1/2})^\nu}{2 \pi i}
  \int_{(a)}  e^s e^{- \pi^2 \ell^2 N / s} s^{- \nu -1}\, \dx s,
\end{equation}
and  hence we have
\begin{equation}
\label{J4-eval} 
 \J_4
=
  \frac{N^{k / 2 + 3 / 4}}{\pi^{k + 1}}
  \sum_{\ell \ge 1}
    \frac{J_{k + 3 / 2} (2 \pi \ell N^{1 / 2})}{\ell^{k + 3 / 2}}.
\end{equation}
The absolute convergence of the series in $ \J_4$  is studied
in \S\ref{sums-abs-conv}.

\subsection{Evaluation of $\J_5$}

With the same substitution used before, since
the double exchange between sums and the line integral
is justified in \S\ref{exchange-double-sum-ell-rho}
for $k>1$, we see that
\begin{equation*}
%\label{J5-def}
  \J_5
  :=
  - \pi^{1 / 2}
  \sum_{\rho}
        \Gamma(\rho)
        N^{k + 1 / 2 + \rho}
        \sum_{\ell \ge 1}
      \Bigl(
        \frac 1{2 \pi i}
        \int_{(1)}
          e^s e^{- \ell^2 \pi^2 N / s} s^{- k - 3 / 2 - \rho} \, \dx s
      \Bigr).
\end{equation*}
Using  \eqref{J-nu}   we get
\begin{equation}
  \J_5
  =
\label{J5-eval}
  -
  \frac{N^{k / 2 + 1 / 4}}{\pi^k}
  \sum_{\rho} \Gamma(\rho) \frac{N^{\rho / 2}}{\pi^\rho}
    \sum_{\ell \ge 1}
      \frac{J_{k + 1 / 2 + \rho} (2 \pi \ell N^{1 / 2})}
           {\ell^{k + 1 / 2 + \rho}}.
\end{equation}
In this case the absolute convergence of the series in $\J_{5}$ 
is more delicate; such a treatment is again described 
in \S\ref{sums-abs-conv}.
 
\medskip
Substituting \eqref{J4-eval}-\eqref{J5-eval} in
\eqref{I2-splitting} we have
\begin{equation}
  \I_2
  =
  \frac{N^{k / 2 + 3 / 4}}{\pi^{k + 1}}
  \sum_{\ell \ge 1}
    \frac{J_{k + 3 / 2} (2 \pi \ell N^{1 / 2})}{\ell^{k + 3 / 2}} 
  -
\label{I2-eval}
  \frac{N^{k / 2 + 1 / 4}}{\pi^k}
  \sum_{\rho}  \frac{\Gamma(\rho)N^{\rho / 2}}{\pi^\rho}
    \sum_{\ell \ge 1}
      \frac{J_{k + 1 / 2 + \rho} (2 \pi \ell N^{1 / 2})}
           {\ell^{k + 1 / 2 + \rho}}.
\end{equation}
Finally, inserting \eqref{final-eval-I1} and \eqref{I2-eval} into
\eqref{hl-splitting} and \eqref{first-step-hl} we finally obtain
\begin{align}
\notag
  \sum_{n \le N} r_{\textit{HL}}(n) \frac{(N - n)^k}{\Gamma(k + 1)}
  &=
  \frac{\pi^{1 / 2}}2 \frac{N^{k + 3 / 2}}{\Gamma(k + 5 / 2)}
  -
  \frac 12 \frac{N^{k + 1}}{\Gamma(k + 2)}
  -
  \frac{\pi^{1 / 2}}2
  \sum_{\rho}
    \frac{\Gamma(\rho)}{\Gamma(k + 3 / 2 + \rho)} N^{k + 1 / 2 + \rho} \\
\notag
  &+
  \frac 12
  \sum_{\rho}
    \frac{\Gamma(\rho)}{\Gamma(k + 1 + \rho)} N^{k + \rho}
  +
  \frac{N^{k / 2 + 3 / 4}}{\pi^{k + 1}}
  \sum_{\ell \ge 1}
    \frac{J_{k + 3 / 2} (2 \pi \ell N^{1 / 2})}{\ell^{k + 3 / 2}} \\
\label{expl-form-HL-bis}
  &-
  \frac{N^{k / 2 + 1 / 4}}{\pi^k}
  \sum_{\rho} \Gamma(\rho) \frac{N^{\rho / 2}}{\pi^\rho}
    \sum_{\ell \ge 1}
      \frac{J_{k + 1 / 2 + \rho} (2 \pi \ell N^{1 / 2})}
           {\ell^{k + 1 / 2 + \rho}}
  +
  \Odip{k}{N^{k}},
\end{align}
for $k > 1$. 
Theorem \ref{CesaroHL-average} follows dividing
\eqref{expl-form-HL-bis} by $N^{k}$.

\section{Lemmas}
We recall some basic facts in complex analysis.
First, if $z = a + i y$ with $a > 0$, we see that for complex $w$ we
have
\begin{align*}
  z^{-w}
  &=
  \vert z \vert^{-w} \exp( - i w \arctan(y / a)) \\
  &=
  \vert z \vert^{-\Re(w) - i \Im(w)} \exp( (- i \Re(w) + \Im(w)) \arctan(y / a))
\end{align*}
so that
\begin{equation}
\label{z^w}
  \vert z^{-w} \vert
  =
  \vert z \vert^{-\Re(w)} \exp(\Im(w) \arctan(y / a)).
\end{equation}
We also recall that, uniformly for $x \in [x_1, x_2]$, with $x_1$ and
$x_2$ fixed, and for $|y| \to +\infty$, by the Stirling formula we have
\begin{equation}
\label{Stirling}
  \vert \Gamma(x + i y) \vert
  \sim
  \sqrt{2 \pi}
  e^{- \pi |y| / 2} |y|^{x - 1 / 2},
\end{equation}
see, e.g., Titchmarsh \cite[\S4.42]{Titchmarsh1988}.

We will need the  following lemmas from Languasco-Zaccagnini \cite{LanguascoZ2012f}.

\begin{Lemma}[See Lemma 1 of \cite{LanguascoZ2012f}] 
\label{Linnik-lemma2}
Let $z = a + iy$, where $a > 0$ and $y \in \R$.
Then
\begin{equation*}
%\label{expl-form-Stilde-strong}
  \widetilde{S}(z)
  =
  \frac{1}{z}
  -
  \sum_{\rho}z^{-\rho} \Gamma(\rho)
  +
  E(a,y)
\end{equation*} 
where $\rho = \beta + i\gamma$ runs over the non-trivial zeros of
$\zeta(s)$ and
\begin{equation}
\label{expl-form-err-term-strong}
  E(a,y)
  \ll
  \vert z \vert^{1/2}
  \begin{cases}
    1 & \text{if $\vert y \vert \leq a$} \\
    1 +\log^2 (\vert y\vert/a) & \text{if $\vert y \vert > a$.}
  \end{cases}
\end{equation}
\end{Lemma}
\begin{Lemma}[See Lemma 2 of \cite{LanguascoZ2012f}] 
\label{series-int-zeros}
Let $\rho=\beta + i \gamma$ run over the non-trivial zeros of the Riemann
zeta-function and $\alpha > 1$ be a parameter.
The series
\[
  \sum_{\rho \colon \gamma > 0}
  \gamma^{\beta-1/2}
    \int_1^{+\infty} \exp\Bigl( - \gamma \arctan\frac 1u \Bigr)
      \frac{\dx u}{u^{\alpha+\beta}}
\]
converges provided that $\alpha > 3/2$.
For $\alpha \le 3/2$ the series does not converge.
The result remains true if we insert in the integral a factor
$(\log u)^c$, for any fixed $c \ge 0$.
\end{Lemma}
%
%\begin{Proof}
%Setting $y = \arctan(1 / u)$, for any real $\gamma > 0$ we have
%%
%\begin{align*}
%  \int_1^{+\infty} \exp\Bigl( -\gamma \arctan\frac 1u \Bigr)
%    \frac{\dx u}{u^{\alpha+\beta}}
%%  &=
%%  \int_0^{\pi / 4}
%%    \exp(-\gamma y) \, \frac{\dx y}{(\cotg y)^\alpha (\sin y)^2} \\
%  &=
%  \int_0^{\pi / 4}
%    \exp(-\gamma y) \,
%    \frac{(\sin y)^{\alpha+\beta - 2}}{(\cos y)^{\alpha+\beta}} \, \dx y\\
%  &\ll_\alpha
%  \int_0^{\pi / 4}
%    \exp(-\gamma y) \, y^{\alpha+\beta - 2} \, \dx y \\
%  &=
%  \gamma^{1 - \alpha - \beta}
%  \int_0^{\pi \gamma / 4}
%    \exp(-w) \, w^{\alpha+\beta - 2} \, \dx w \\
%  &\le
%  \gamma^{1 - \alpha - \beta} \,
%  \Bigl(\Gamma(\alpha-1)+\Gamma(\alpha)\Bigr),
%\end{align*}
%%
%since $0 < \beta < 1$.
%This shows that the series over $\gamma$ converges for $\alpha > 3/2$.
%For $\alpha = 3/2$ essentially the same computation shows that the
%integral is $\gg \gamma^{-1/2 - \beta}$ and it is well known that in
%this case the series over zeros diverges.
%The other assertions are proved in the same way.
%\end{Proof}
%
\begin{Lemma}[See Lemma 3 of \cite{LanguascoZ2012f}] 
\label{series-int-zeros-alt-sign}
Let $\alpha > 1$, $z=a+iy$, $a\in(0,1)$ and $y\in \R$.
Let further $\rho=\beta+i\gamma$ run over the non-trivial zeros of
the Riemann zeta-function.  We have
\[
  \sum_{\rho}
    \vert \gamma\vert ^{\beta-1/2}
    \int_{\Y_1 \cup \Y_2} \exp\Bigl(\gamma \arctan\frac{y}{a} - \frac\pi2 \vert \gamma \vert\Bigr)
      \frac{\dx y}{\vert z \vert ^{\alpha+\beta}}
  \ll_{\alpha}
  a^{-\alpha},
\]
where $\Y_1=\{y\in \R\colon y\gamma \leq 0\}$ and 
$\Y_2=\{y\in [-a,a] \colon y\gamma > 0\}$.
The result remains true if we insert in the integral a factor
$(\log (\vert y\vert /a))^c$, for any fixed $c \ge 0$.
\end{Lemma}
%
%\begin{Proof}
%We first work on $\Y_1$.
%By symmetry, we may assume that $\gamma > 0$.
%For $y \in(-\infty, 0]$ we have
%$\gamma \arctan(y/a) -\frac \pi2 \vert \gamma\vert  \le - \frac \pi2 \vert \gamma\vert $
%and hence the quantity we are estimating becomes
%\[
%  \sum_{\rho \colon \gamma > 0} 
%  \gamma^{\beta-1/2}
%  \exp\Bigl( -\frac \pi2 \gamma \Bigr)
%    \int_{-\infty}^0 \frac{\dx y}{\vert z \vert ^{\alpha+\beta}}
%  \ll_\alpha
%  \sum_{\rho \colon \gamma > 0} 
%  \gamma^{\beta-1/2}
%  \exp\Bigl( -\frac \pi2 \gamma \Bigr)
%  a^{1-\alpha-\beta}
%  \ll_{\alpha}
%  a^{-\alpha},
%\]
%using $0<\beta<1$, standard zero-density estimates and \eqref{z^-1}.
%We consider now the integral over $\Y_2$.
%Again by symmetry we can assume that $\gamma > 0$ and so we get
%%
%\begin{align*}
%  \sum_{\rho \colon \gamma > 0}
%   \gamma ^{\beta-1/2}
%    \int_0^{a} \exp\Bigl( \gamma (\arctan \frac{y}{a} -\frac \pi2) \Bigr)
%     \frac{\dx y}{\vert z \vert ^{\alpha+\beta}}
%  &\ll
%  \sum_{\rho \colon \gamma > 0}
%  \gamma ^{\beta-1/2}
%  \exp\Bigl( -\frac \pi4 \gamma \Bigr)
%  \int_0^{a} \frac{\dx y}{\vert z \vert ^{\alpha+\beta}} \\
%  &\ll_\alpha
%  \sum_{\rho \colon \gamma > 0}
%  \gamma ^{\beta-1/2} \exp\Bigl(  -\frac \pi4 \gamma \Bigr)
%  a^{1-\alpha-\beta} 
%  \ll_\alpha
%  a^{-\alpha}
%\end{align*}
%%
%arguing as above.
%The other assertions are proved in the same way.
%\end{Proof}

\section{Interchange of the series over zeros with the line integral} 
\label{exchange-rho-integral}

We need $k>0$ in this section. We have to establish the convergence of
\begin{equation}
\label{conv-integral-J2-J3}
  \sum_{\rho}
    \vert \Gamma(\rho)\vert  \int_{(1/N)} \vert e^{N z}\vert  \vert z\vert ^{- k - 3/2} \vert z^{- \rho}\vert  \, \vert \dx z\vert,
\end{equation}
where, as usual, $\rho=\beta + i \gamma$ runs over the non-trivial 
zeros of the Riemann zeta-function.
By \eqref{z^w} and the Stirling formula \eqref{Stirling}, we are left
with estimating
\begin{equation}
\label{conv-integral-J2-J3-1}
  \sum_{\rho}
    \vert \gamma\vert ^{\beta - 1/2}
    \int_{\R} \exp\Bigl( \gamma \arctan(N y) -\frac \pi2 \vert \gamma\vert \Bigr)
             \frac{\dx y}{\vert z \vert ^{k + 3/2 +\beta}}. 
\end{equation}
We have just to consider the case $\gamma y >0$, $\vert y \vert > 1/N$
since in the other cases the total contribution is  $\ll_k N^{k + 1 }$
by Lemma \ref{series-int-zeros-alt-sign} with $\alpha=k+3/2$ and $a=1/N$.
By symmetry, we may assume that $\gamma > 0$.
We have that the integral in \eqref{conv-integral-J2-J3-1} is
\begin{align*}
    &
  \ll
  \sum_{\rho \colon \gamma > 0}
    \gamma ^{\beta - 1/2}
    \int_{1 / N}^{+\infty} \exp\Bigl( - \gamma \arctan\frac 1{N y} \Bigr)
      \frac{\dx y}{y^{k + 3/2 +\beta}} \\
  &=
  N^{k}
  \sum_{\rho \colon \gamma > 0}
    N^{\beta}
    \gamma ^{\beta - 1/2}
    \int_1^{+\infty} \exp\Bigl( - \gamma \arctan\frac 1u \Bigr)
      \frac{\dx u}{u^{k + 3/2 +\beta}}.
\end{align*}
For $k > 0 $ this is $\ll_k N^{k + 1}$ by
Lemma~\ref{series-int-zeros}. This implies
that the integrals in \eqref{conv-integral-J2-J3-1} and in \eqref{conv-integral-J2-J3} 
are both  $\ll_k N^{k + 1}$ and
hence this exchange step is fully justified.

\section{Interchange of the series over $\ell$ with the line integral}
\label{exchange-ell-integral}

We need $k> - 1/2$ in this section.
We  have to  establish the convergence of
\begin{equation}
\label{conv-integral-J4}
  \sum_{\ell \ge 1}
    \int_{(1/N)} 
    \vert e^{N z}\vert  \vert z\vert ^{- k - 5/2}
    e^{-\pi^2 \ell^2 \Re(1/z)}  \, \vert \dx z\vert .
\end{equation}
A trivial computation gives 
\begin{equation}
\label{real-part-estim}
\Re(1/z)
=
\frac{N}{1+N^2y^2}
\gg
\begin{cases}
N & \text{if}\ \vert y \vert \leq 1/N \\
1/(Ny^2) & \text{if}\ \vert y \vert > 1/N.
\end{cases}
\end{equation}
By \eqref{real-part-estim}, we can write that
the  quantity in \eqref{conv-integral-J4} is
\begin{equation} 
    \label{J4-split}
\ll 
  \sum_{\ell \ge 1}
    \int_{0}^{1/N}  
    \frac{e^{- \ell^2 N}}{\vert z\vert ^{k + 5/2}}  \, \dx y
    +
 \sum_{\ell \ge 1}
    \int_{1/N}^{+\infty} 
    \frac{e^{- \ell^2 / (Ny^2)}}{\vert z\vert ^{k + 5/2}}   \, \dx y
    =
U_1+U_2,
\end{equation}
say, since the $\pi^2$ factor in the exponential function
is negligible. Using \eqref{omega-estim}-\eqref{z^-1},  we have 
\begin{equation}
\label{U1-estim}
U_1
\ll
 N^{k+3/2}
\omega_{2} (N) 
 \ll  
 N^{k+1}
\end{equation}
and
\begin{align}
\notag
U_2&
\ll
\sum_{\ell \ge 1}
    \int_{1/N}^{+\infty} 
    \frac{e^{- \ell^2 / (Ny^2)}}{y^{k + 5/2}}  \, \dx y
  \ll
N^{k/2+3/4} \sum_{\ell \ge 1}
\frac{1}{\ell^{k+3/2}}
    \int_{0}^{\ell^2 N} u^{ k/2 - 1/4} 
    e^{- u}  \, \dx u 
\\
\label{U2-estim}
&\leq
\Gamma \Bigl( \frac{2k+3}{4} \Bigr)
N^{k/2+3/4}  \sum_{\ell \ge 1}
\frac{1}{\ell^{k+3/2}}
\ll_{k}
N^{k/2+3/4}
\end{align}
provided that $k>-1/2$,
where  we used the substitution
$u=\ell^2 / (Ny^2)$.

Inserting \eqref{U1-estim}-\eqref{U2-estim}
into \eqref{J4-split} we get, for $k> - 1/2$, that
the quantity in \eqref{conv-integral-J4}
is $\ll N^{k+1}$.

\section{Interchange of the double series over zeros with the line integral} 
\label{exchange-double-sum-ell-rho}

We need $k>1$ in this section.
We first have to  establish the convergence of
\begin{equation}
\label{conv-integral-J5}
  \sum_{\ell \ge 1}
    \int_{(1/N)} 
         \vert \sum_{\rho} \Gamma(\rho)z^{- \rho} \vert
    \vert e^{N z}\vert  \vert z\vert ^{- k - 3/2}
    e^{-\pi^2 \ell^2 \Re(1/z)}  \, \vert \dx z\vert .
\end{equation}
Using the PNT and \eqref{expl-form-err-term-strong}, we first remark
that
\begin{align}
\notag
  \Bigl\vert \sum_{\rho} z^{-\rho} \Gamma(\rho) \Bigr\vert 
  &=
  \Bigl\vert  \Stilde(z) - \frac 1z -E(y,\frac{1}{N}) 
  \Bigr\vert 
  \ll
  N  + \frac{1}{\vert z\vert}
  +
  \Bigl\vert   E(y,\frac{1}{N}) \Bigr\vert
\\
\label{sum-over-rho-new}
  &\ll
  \begin{cases}
    N & \text{if $\vert y \vert \leq 1/N$,} \\
    \vert z\vert ^{-1} + \vert z\vert ^{1/2}
    \log^2 (2N\vert y\vert) & \text{if $\vert y \vert > 1/N$.}
  \end{cases}
\end{align}

By \eqref{real-part-estim} and \eqref{sum-over-rho-new}, we can write that
the  quantity in \eqref{conv-integral-J5} is
\begin{align}
\notag 
&\ll N
  \sum_{\ell \ge 1}
    \int_{0}^{1/N}  
    \frac{e^{- \ell^2 N}}{\vert z\vert ^{k + 3/2}}  \, \dx y
    &+
 \sum_{\ell \ge 1}
    \int_{1/N}^{+\infty}   
    \frac{e^{- \ell^2 / (Ny^2)} }{\vert z\vert ^{k + 5/2}}  \, \dx y
    +
    \sum_{\ell \ge 1}
    \int_{1/N}^{+\infty}   \log^{2}(2Ny)
    \frac{e^{- \ell^2 / (Ny^2)}}{\vert z\vert ^{k + 1}}   \, \dx y
    \\&
    \label{HL-split}
    =
V_1+V_2+V_{3},
\end{align}
say.
$V_{1}$ and $V_{2}$ can be estimated exactly as $U_{1}, U_{2}$
in \S\ref{exchange-ell-integral}; hence we have 
\begin{equation}
\label{V1-V2-estim}
V_{1} + V_{2} \ll_{k} N^{k+1}
\end{equation}
provided that  $k> - 1/2$.

Using the substitution $u=\ell^2 / (Ny^2)$,
we obtain
\[
V_3
\ll
 \sum_{\ell \ge 1}
    \int_{1/N}^{+\infty} \log^{2}(2Ny) 
    \frac{e^{- \ell^2 / (Ny^2)}}{y^{k + 1}}   \, \dx y
=
\frac{N^{k/2}}{8} \sum_{\ell \ge 1}
\frac{1}{\ell^k}
    \int_{0}^{\ell^2 N} u^{ k/2 - 1} 
    \log^{2}\Bigl( \frac{4\ell^{2}N}{u} \Bigr)
    e^{- u}  \, \dx u.
\]
Hence a direct computation shows that
\begin{align}
\notag
V_3
&\ll
N^{k/2} \sum_{\ell \ge 1}
\frac{\log^{2}(\ell N)}{\ell^k}
    \int_{0}^{\ell^2 N} u^{ k/2 - 1}     e^{- u}  \, \dx u
    +
N^{k/2} \sum_{\ell \ge 1}
\frac{1}{\ell^k}
     \int_{0}^{\ell^2 N} u^{ k/2 - 1}  \log^{2} (u)\,
    e^{- u}  \, \dx u
\\
\label{V3-estim}
&\ll_{k}
\Gamma(k/2)
N^{k/2} \sum_{\ell \ge 1}
\frac{\log^{2}(\ell N)}{\ell^k}
    +
N^{k/2}
\ll_{k}
N^{k/2}\log^{2} N
\end{align}
provided that $k>1$.
Inserting \eqref{V1-V2-estim}-\eqref{V3-estim}
into \eqref{HL-split} we get, for $k>1$, that
the quantity in \eqref{conv-integral-J5}
is $\ll N^{k+1}$.

Now we have to  establish the convergence of
\begin{equation}
\label{conv-integral-5}
  \sum_{\ell \ge 1}
     \sum_{\rho}\vert \Gamma(\rho) \vert 
    \int_{(1/N)} 
   \vert e^{N z}\vert  \vert z\vert ^{- k - 3/2}
         \vert z^{- \rho} \vert
    e^{-\pi^2 \ell^2 \Re(1/z)}  \, \vert \dx z\vert .
\end{equation}

By symmetry, we may assume that $\gamma > 0$.
For $y \in(-\infty, 0]$ we have
$\gamma \arctan(y/a) -\frac \pi2  \gamma  \le - \frac \pi2  \gamma$.
Using \eqref{real-part-estim}, \eqref{z^-1} and the Stirling formula \eqref{Stirling}, 
the quantity we are estimating becomes
\begin{align}
\notag
&\ll
  \sum_{\ell \ge 1}
  \sum_{\rho \colon \gamma>0}
  \gamma^{\beta-1/2}
   \exp\Bigl(  -\frac \pi2 \gamma\Bigr)
   \Bigl(
              \int_{-1/N}^{0}  
         N^{k + 3/2 + \beta}\ e^{- \ell^2 N}  \, \dx y
        +
        \int_{-\infty}^{-1/N}
         \frac{e^{- \ell^2 / (Ny^2)}}{\vert y \vert^{k + 3/2 + \beta}}  \, \dx y
        \Bigr)
    \\
    \notag
   & \ll_{k}
N^{k+3/2}
 \sum_{\ell \ge 1}    e^{- \ell^2 N} 
\sum_{\rho \colon \gamma>0}
\gamma^{\beta-1/2}
\exp\Bigl(  -\frac \pi2 \gamma\Bigr)
\\
\notag
&\hskip1cm
+
N^{k/2+1/4} \sum_{\ell \ge 1}
\frac{1}{\ell^{k+1/2}}
\sum_{\rho \colon \gamma>0}
\frac{N^{\beta/2}}{\ell^{\beta}}
\gamma^{\beta-1/2}
\exp\Bigl(  -\frac \pi2 \gamma\Bigr)
    \int_{0}^{\ell^2 N} u^{ k/2 - 3/4 + \beta/2} 
    e^{- u}  \, \dx u 
    \\
    \notag
    &
    \ll_{k}
     N^{k+1}
     + 
\Bigl(
  \max_{\beta} \Gamma\Bigl( \frac{\beta}{2} +\frac{k}2+ \frac14\Bigr) 
\Bigr) 
N^{k/2+3/4} \sum_{\ell \ge 1}
\frac{1}{\ell^{k+1/2}}
\sum_{\rho \colon \gamma>0} 
\gamma^{\beta-1/2}
\exp\Bigl(  -\frac \pi2 \gamma\Bigr)
    \\
    \label{easy-case}
    & \ll_k  
 N^{k+1}
\end{align}
provided that $k>1/2$, 
where  we used the substitution $u = - \ell^2 / (Ny^2)$, 
 \eqref{omega-estim} and standard density estimates.

Let now $y>0$.
Using the Stirling formula \eqref{Stirling} and \eqref{real-part-estim} we can write that
the  quantity in \eqref{conv-integral-5} is
\begin{align}
\notag
\ll
  \sum_{\ell \ge 1}
  \sum_{\rho \colon \gamma>0}
 & \gamma^{\beta-1/2}
\exp\Bigl(  -\frac \pi4 \gamma\Bigr)
              \int_{0}^{1/N}  
        \frac{e^{- \ell^2 N}}{\vert z\vert ^{k + 3/2 + \beta}}  \, \dx y
   \\ 
   \notag
   &+
\sum_{\ell \ge 1}
  \sum_{\rho \colon \gamma>0}
  \gamma^{\beta-1/2}
      \int_{1/N}^{+\infty}  
\exp\Bigl( \gamma (\arctan(N y) -\frac \pi2) \Bigr)  
    \frac{e^{- \ell^2 / (Ny^2)}}{\vert z\vert ^{k + 3/2 + \beta}}  \, \dx y
    \\
    \label{HL-split-1}
    &
    =
W_1+W_2,
\end{align}
say.
Using \eqref{z^-1} and \eqref{omega-estim},  we have that
\begin{equation}
\label{W1-estim}
W_1
\ll
N^{k+3/2}
 \sum_{\ell \ge 1}    e^{- \ell^2 N} 
\sum_{\rho \colon \gamma>0}
\gamma^{\beta-1/2}
\exp\Bigl(  -\frac \pi4 \gamma\Bigr)
 \ll  
 N^{k+1}
\end{equation}
by standard density estimates.
Moreover we get
\begin{align*}
W_2
&\ll
\sum_{\ell \ge 1}
\sum_{\rho \colon \gamma>0}
\gamma^{\beta-1/2}
        \int_{1/N}^{+\infty} y^{- k - 3/2 - \beta}  
\exp\Bigl( - \frac{\gamma}{N y} - \frac{\ell^2}{Ny^2} \Bigr)  \, \dx y
\\
&=
N^{k/2+1/4} 
\sum_{\ell \ge 1}
\frac{1}{\ell^{k+1/2}}
\sum_{\rho \colon \gamma>0}
\frac{N^{\beta/2} \gamma^{\beta-1/2}} {\ell^{\beta}}
    \int_{0}^{\ell \sqrt{N}} v^{ k - 1/2 + \beta} 
    \exp\Bigl(-\frac{\gamma v}{\ell\sqrt{N}}- v^2\Bigr) 
    \, \dx v,
\end{align*}
in which we used the substitution
$v^2=\ell^2 / (Ny^2)$.
Remark now that, for $k>1$, 
we can set $\eps=\eps(k)=(k-1)/2>0$ and that 
$k -\eps =(k+1)/2>1$.
We further remark that $\max_{v} (v^{k -\eps}e^{-v^2})$
is attained at $v_0=((k -\eps)/2)^{1/2}$, and hence we obtain,
for $N$ sufficiently large, that
\[
W_2
\ll_{k}
N^{k/2+1/4}
  \sum_{\ell \ge 1}
\frac{1}{\ell^{k+1/2}}
\sum_{\rho \colon \gamma>0}
\frac{N^{\beta/2} \gamma^{\beta-1/2}} {\ell^{\beta}}
    \int_{0}^{\ell \sqrt{N}}  
    v^{\beta-1/2+\eps}\exp\Bigl(-\frac{\gamma v}{\ell\sqrt{N}}\Bigr)
     \, \dx v .
\]
Making the substitution $u= \gamma v/(\ell\sqrt{N})$  we have
\begin{align}
\notag
W_2
&\ll_{k}
N^{k/2+1/2+\eps/2}
  \sum_{\ell \ge 1}
\frac{1}{\ell^{k -\eps}}
\sum_{\rho \colon \gamma>0}
\frac{N^{\beta}} {\gamma^{1+\eps}}
    \int_{0}^{\gamma}  
    u^{\beta-1/2+\eps} e^{- u}
     \, \dx u\\
\label{W2-estim}
&\ll_{k}
N^{k/2+3/2+\eps/2}
  \sum_{\ell \ge 1}
\frac{1}{\ell^{k -\eps}}
\sum_{\rho \colon \gamma>0}
\frac{ 1 } {\gamma^{1+\eps}}
\Bigl(
\max_{\beta} \Gamma \Bigl( \beta + \frac12+\eps \Bigr)  
\Bigr)  
\ll_{k}
N^{k/2+3/2+\eps/2},
\end{align}
by standard density estimates. 

Inserting \eqref{W1-estim}-\eqref{W2-estim}
into \eqref{HL-split-1} and recalling \eqref{easy-case},
we get, for $k>1$, that
the quantity in \eqref{conv-integral-5}
is $\ll  N^{k+1}$.

\section{Absolute convergence of $\J_{4}$ and $\J_{5}$}
\label{sums-abs-conv}

\medskip
To study the absolute convergence of the series in $ \J_4$ we first remark that,
by \eqref{Bessel-def} and \eqref{J4-eval}, we get
\[
  \sum_{\ell \ge 1}
    \frac{\vert J_{k + 3 / 2} (2 \pi \ell N^{1 / 2}) \vert }{\ell^{k + 3 / 2}}
    \ll_{k}
    N^{- k / 2 - 3 / 4}
    \sum_{\ell \ge 1}
    \int_{(1/N)} 
    \vert e^{N z}\vert  \vert z\vert ^{- k - 5/2}
    e^{-\pi^2 \ell^2 \Re(1/z)}  \, \vert \dx z\vert 
\]
which is the quantity in \eqref{conv-integral-J4}.
So the argument in \S\ref{exchange-ell-integral}
also proves that the series in $\J_{4}$
converges absolutely for  $k> -1/2$.

In fact a   more direct argument leads to a better estimate on $k$.
Using, for $\nu>0$ fixed, $u\in \R$ and $u\to+\infty$, the estimate
\begin{equation*}
%\label{Berndt-estimate}
\vert J_\nu(u) \vert
\ll_\nu 
 u^{-1/2}
\end{equation*}
which immediately follows from \eqref{Lebedev-asymp} 
(or from eq.~(2.4) of Berndt \cite{Berndt1975}), and performing
 a direct computation, we obtain that $ \J_4$
converges absolutely for $k> -1$ (and for $N$ sufficiently large)
and that $\J_4\ll_{k} N^{(k+1)/2}$.

For the study of the absolute convergence of the series in $\J_5$ we have a different situation. In this case the direct argument needs a more careful estimate of the Bessel functions involved since both $\nu$ and $u$ are not fixed and, in fact, unbounded. 
In fact it is easy to see that \eqref{Lebedev-asymp} can be used only if 
$\nu \in \C$ is bounded, but we are not in this case
since $\nu=k + 1 / 2 + \rho$, where $\rho$ is a nontrivial zero of the 
Riemann $\zeta$-function. On the other hand, \eqref{Poisson-int-rep} 
can be used only for $u$ bounded, but again this is not our case
since $u= 2 \pi \ell N^{1 / 2}$ and $\ell$ runs up to infinity.
Moreover, the use of the asymptotic relations for $J_\nu(u)$
when  $\nu\in\C$ and $u\in\R$ are both ``large'' seems to be very complicated
in this setting. 

So it turned out that the best direct approach we are able to perform is the following.
By a double partial integration on \eqref{Poisson-int-rep}, we immediately get 
\begin{align} 
\notag
J_\nu(u)
&=
\frac{2(u/2)^{\nu}(2\nu-1)}{\pi^{1/2}u^{2}\Gamma(\nu+1/2)}
\int_{0}^{1} 
\Bigl(
1 - \frac{(2\nu-3)t^{2}}{1-t^2}
\Bigr)
(1-t^2)^{\nu-3/2} 
\cos (ut)\ \dx t
\\
\notag
&\ll_{\Re(\nu)}
\frac{\vert u\vert ^{\Re(\nu)-2}\vert 2\nu-1\vert}{\vert \Gamma(\nu+1/2) \vert}
\int_{0}^{1} 
\Bigl( 1 + \vert  2\nu-3 \vert  
\Bigr)
\vert \cos (ut) \vert \ \dx t
\\&
\label{Poisson-estim2}
\ll_{\Re(\nu)}
\frac{\vert \nu \vert^{2}\vert u\vert ^{\Re(\nu)-2}}{\vert  \Gamma(\nu+1/2) \vert},
\end{align}
where the last two estimates hold for $\Re(\nu)>3/2$
and $u>0$. 
Inserting \eqref{Poisson-estim2} into \eqref{J5-eval} and using 
 the Stirling formula \eqref{Stirling}, 
a direct computation shows  the absolute convergence of the double sum 
in $\J_{5}$  for $k>2$ (and for $N$ sufficiently large).

Unfortunately, such a condition on $k$ is worse than the one we have in 
\S\ref{exchange-double-sum-ell-rho}. So, coming back to the Sonine 
representation of the Bessel functions \eqref{Bessel-def} on the line
$\Re(s)=1$ and using the usual substitution $s=Nz$, 
to study the absolute convergence of the double sum in $\J_{5}$
we are led to consider the quantity
\begin{align*}
\sum_{\rho} 
 \Bigl\vert
  \Gamma(\rho) \frac{N^{\rho / 2}}{\pi^\rho} 
  \Bigr\vert
  &  \sum_{\ell \ge 1}
       \Bigl\vert
       \frac{ J_{k + 1 / 2 + \rho} (2 \pi \ell N^{1 / 2}) }
           {\ell^{k + 1 / 2 + \rho}} 
           \Bigr\vert
           \\
&           
    \ll_{k}
    N^{-k/2-1/4}
     \sum_{\rho}\vert \Gamma(\rho) \vert 
     \sum_{\ell \ge 1}  
    \int_{(1/N)} 
   \vert e^{N z}\vert  \vert z\vert ^{- k - 3/2}
         \vert z^{- \rho} \vert
    e^{-\pi^2 \ell^2 \Re(1/z)}  \, \vert \dx z\vert,
\end{align*}
which is very similar to the one in \eqref{conv-integral-5}
(the sums are interchanged). It is not hard to see that 
the argument used in \eqref{conv-integral-5}-\eqref{W2-estim} can be applied in this case
too. It shows that the double series in $\J_{5}$ converges absolutely for 
$k>1$ and this condition fits now with the one we have in  
\S\ref{exchange-double-sum-ell-rho}.

\newpage
%%\bibliographystyle{amsplain}
%\bibliographystyle{amsplain-nobysame}
%\bibliography{teonum}

\providecommand{\bysame}{\leavevmode\hbox to3em{\hrulefill}\thinspace}
\providecommand{\MR}{\relax\ifhmode\unskip\space\fi MR }
\providecommand{\MRhref}[2]{%
  \href{http://www.ams.org/mathscinet-getitem?mr=#1}{#2}
}
\providecommand{\href}[2]{#2}

\vskip0.5cm
%\bigskip
\noindent
\begin{tabular}{l@{\hskip 24mm}l}
Alessandro Languasco               & Alessandro Zaccagnini\\
Universit\`a di Padova     & Universit\`a di Parma\\
Dipartimento di Matematica & Dipartimento di Matematica \\
Via Trieste 63                & Parco Area delle Scienze, 53/a \\
35121 Padova, Italy            & 43124 Parma, Italy\\
{\it e-mail}: languasco@math.unipd.it        & {\it e-mail}:
alessandro.zaccagnini@unipr.it  
\end{tabular}

\end{document}